\documentclass[11pt,a4paper,final]{amsart}
\usepackage{amsfonts,amsmath,amssymb}
\usepackage{pstricks,pst-plot}
\numberwithin{equation}{section}

\renewcommand{\epsilon}{\varepsilon}

\DeclareSymbolFont{SY}{U}{psy}{m}{n}
\DeclareMathSymbol{\emptyset}{\mathord}{SY}{'306}

\DeclareMathOperator{\Ran}{Ran} 
 \DeclareMathOperator{\Dom}{Dom}
\DeclareMathOperator{\sign}{sign}

\DeclareMathOperator{\spec}{spec}

\DeclareMathSymbol{\newtimes}{\mathbin}{SY}{'264}

\newcommand{\R}{\mathbb{R}}

\newcommand{\C}{\mathbb{C}}
\newcommand{\Z}{\mathbb{Z}}
\newcommand{\N}{\mathbb{N}}

\newcommand{\fH}{\mathfrak{H}}

\newcommand{\fL}{\mathfrak{L}}

\newcommand{\fK}{\mathfrak{K}}

\newcommand{\fa}{\mathfrak{a}}
\newcommand{\fb}{\mathfrak{b}}
\newcommand{\fh}{\mathfrak{h}}

\newcommand{\fv}{\mathfrak{v}}

\newcommand{\fx}{\mathfrak{x}}
\newcommand{\fy}{\mathfrak{y}}

\newcommand{\cB}{{\mathcal B}}

\newcommand{\cD}{{\mathcal D}}

\newcommand{\cH}{{\mathcal H}}

\newcommand{\ii}{\mathrm{i}}

\newtheorem{theorem}{Theorem}[section]{\bf}{\it}
\newtheorem{hypothesis}[theorem]{Hypothesis}{\bf}{\it}
\newtheorem{proposition}[theorem]{Proposition}{\bf}{\it}
{\bf}{\it}
\newtheorem{example}[theorem]{Example}{\it}{\rm}
\newtheorem{lemma}[theorem]{Lemma}{\bf}{\it}
\newtheorem{remark}[theorem]{Remark}{\it}{\rm}
{\bf}{\it}
{\bf}{\it}
{\bf}{\it}
{\bf}{\it}

\title[Representation Theorems]{Representation Theorems\\ for Indefinite Quadratic Forms Revisited}

\author[L. Grubi\v{s}i\'c]{Luka Grubi\v{s}i\'c}
\address{L.~Grubi\v{s}i\'c,
Department of Mathematics, University of Zagreb, Bijeni\v{c}ka 30,
10000 Zagreb, Croatia}
\email{luka.grubisic@math.hr}

\author[V. Kostrykin]{Vadim Kostrykin}
\address{V.~Kostrykin, FB 08 - Institut f\"{u}r Mathematik,
Johannes Gutenberg-Universit\"{a}t Mainz,
Staudinger Weg 9,
D-55099 Mainz,
Germany}
\email{kostrykin@mathematik.uni-mainz.de}

\author[K. A. Makarov]{Konstantin A.~Makarov}
\address{K.~A.~Makarov, Department of Mathematics, University of
Missouri, Co\-lum\-bia, MO 65211, USA}
\email{makarovk@missouri.edu}

\author[K. Veseli\'c]{Kre\v{s}imir Veseli\'c}
\address{K.~Veseli\'c,
Fakult\"{a}t f\"{u}r Mathematik und Informatik, Fernuniversit\"{a}t Hagen, Postfach 940,
D-58084 Hagen, Germany}
\email{kresimir.veselic@fernuni-hagen.de}



\begin{document}

\maketitle

\begin{abstract}
The first and second representation theorems for sign-indefinite,  not necessarily semi-bounded
quadratic forms are revisited. New straightforward proofs of these theorems are given. A number
of necessary and sufficient conditions ensuring the second representation theorem to hold
are obtained. A new simple and explicit example of a self-adjoint operator for which the second
representation theorem fails to hold is also provided.
\end{abstract}

\footnotetext{This work is supported in part by the Deutsche Forschungsgemeinschaft,
Grant  KO 2936/3-1}

%%%%%%%%%%%%%%%%%%%%%%%%%%%%%%%%%%%%%%%%%%%%%
\section{Introduction}\label{sec:intro}
%%%%%%%%%%%%%%%%%%%%%%%%%%%%%%%%%%%%%%%%%%%%%

In this work we revisit the representation theorems for sign-indefinite,
not necessarily semi\-boun\-ded symmetric sesquilinear forms. Let $\fH$ be a complex Hilbert space with the inner product $\langle\cdot,\cdot\rangle$. We will be dealing with the class of forms given by
\begin{equation}\label{forma}
\fb[x,y]=\langle A^{1/2}x,HA^{1/2}y\rangle,\quad x,y\in \Dom[\fb] = \Dom  (A^{1/2}),
\end{equation}
where $A$ is uniformly positive self-adjoint operator in the Hilbert space $\fH$,
and $H$ is a bounded, not necessarily positive, self-adjoint operator in $\fH$. In perturbation theory, such forms
arise naturally, when the initial (in general sign-indefinite)
form $\fa$ given by
$$\fa[x,y]=\langle A^{1/2}x, J_A A^{1/2}y\rangle,\quad \quad x,y\in\Dom(A^{1/2}),$$
with
$J_A$ a self-adjoint involution commuting with $A$, is perturbed by a form $\fv$ satisfying the
upper bound
\begin{equation}\label{bnd}
 |\fv[x,y]| \leq \beta |\langle A^{1/2}x, A^{1/2}y\rangle|,\quad x,y\in\Dom(A^{1/2}),
\end{equation}
for some $\beta>0$. So that the sesquilinear form
\begin{equation*}
 \fb[x,y]:=\fa[x,y]+ \fv[x,y]
\end{equation*}
can be transformed into the expression given by \eqref{forma} for some
self-adjoint bounded operator $H$. In particular, Dirac-Coulomb operators
fit into this scheme \cite{Nenciu}. In this connection it is also worth mentioning
an alternative approach
to Dirac-like operators developed recently by Esteban and Loss in \cite{Esteban:Loss}.

In this setting, in the framework of a unified approach, we provide new straightforward proofs of the following two assertions (Theorems \ref{thm:main} and \ref{thm:main:2}, respectively):

\medskip

{\it
\begin{itemize}
\item[(i)]
If $H$ has a bounded inverse, then there is a unique self-adjoint boundedly invertible operator $B$ with $\Dom(B)\subset\Dom[\fb]$ \emph{associated} with the form $\fb$,
that is,
\begin{equation}\label{eq:1RT}
\fb[x,y]= \langle x, By\rangle\quad\text{for all}\quad x\in\Dom[\fb],\quad y\in\Dom(B)\subset\Dom[\fb].
\end{equation}

\medskip

\item[(ii)]
If, in addition, the domains of $|B|^{1/2}$ and $A^{1/2}$ agree, that is,
\begin{equation}\label{ravny}
\Dom(|B|^{1/2}) = \Dom(A^{1/2}),
\end{equation}
then the form  $\fb$  is \emph{represented} by $B$ in the sense that
\begin{equation}\label{eq:2RT}
\fb[x,y]= \langle |B|^{1/2} x,
\sign(B)|B|^{1/2} y\rangle\quad\text{for all}
\quad x,y\in\Dom[\fb]=\Dom(|B|^{1/2}),
\end{equation}
with $\sign(B)$ the sign of the operator $B$.
\end{itemize}
}

Representations \eqref{eq:1RT} and \eqref{eq:2RT} are
natural generalizations to the case of indefinite forms
of the First and Second Representation Theorems for semi-bounded sesquilinear forms,
Theorems VI.2.1 and VI.2.23 in \cite{Kato}, respectively.
We remark that the existence and uniqueness of the pseudo-Friedrichs extension
for symmetric operators \cite[Section VI.4]{Kato},
\cite[Section IV.4]{Evans} is a particular case of this result.

A proof of the First Representation Theorem (i) for
indefinite forms can be found in \cite[Theorem 2.1]{Nenciu}.
Results equivalent to \eqref{eq:1RT} have been obtained by McIntosh in \cite{McIntosh:1},
\cite{McIntosh:2} where, in particular, the notion of closedness of a semibounded form to the case of
indefinite forms has been extended.  It is worth mentioning that the
form $\fb$ given by \eqref{forma} is closed in that sense (see Remark \ref{rem:mcintosh}).

The Second Representation Theorem (ii) for indefinite forms, is originally
due to McIntosh  \cite{McIntosh:1}, \cite{McIntosh:2}. He also established that the form-domain stability
criterion \eqref{ravny} is equivalent to the requirement
that the self-adjoint involution $\sign(B)$ leaves $\Dom (A^{1/2})$ invariant. We remark that if $B$ is a semi-bounded operator, the condition \eqref{ravny} holds automatically (cf.~\cite[Theorem VI.2.23]{Kato}).

New proofs of the Representation Theorems (i) and (ii) given in the present work
are straightforward and based on functional-analytic ideas similar
to those used in the familiar proofs of the Representation Theorems
in the semi-bounded case (cf.~\cite[Section VI.2]{Kato}).
Related results, in particular those concerned with the so-called
quasi-definite matrices and operators are contained in \cite{Veselic:2000} and,
quite recently, appeared in \cite{Veselic}.

Our new results related to the Representation Theorems (i) and (ii) are as follows.

As a consequence of (i), we prove the First Representation Theorem for block operator
matrices defined as quadratic forms, provided that the diagonal part of the matrix has a
bounded inverse and that the off-diagonal form perturbation is relatively bounded with
respect to a closed positive definite form generated
by the diagonal entries of the matrix (Theorem \ref{off} below).
This theorem provides a far reaching
generalization of a result obtained previously by Konstantinov and Mennicken in \cite{KonMenn}.

In this context, we also revisit the Lax-Milgram theory for
coercive closed forms (cf., \cite[Theorem IV.1.2]{Evans}) and show that the coercivity hypothesis yields the
representation \eqref{forma} (Theorem \ref{cor:1} below).

With regard to the Second Representation Theorem (ii), we obtain a number of new necessary and
sufficient conditions for the domains $\Dom(|B|^{1/2})$ and $\Dom(A^{1/2})$ to coincide
(Theorem \ref{lem:n} below). Answering a question raised by A.~McIntosh in  \cite{McIntosh:2},
we also provide an explicit example (Example \ref{exam}) of a form $\fb$, that is $0$-closed in the sense of McIntosh (see Remark \ref{rem:mcintosh}),
but not represented by its associated operator $B$. Consequently, the Second
Representation Theorem fails to hold if the condition \eqref{ravny} is violated.
In particular, we show that the  $A$-form boundedness of the operator $B$
does not yield the $A$-form boundedness for its absolute value $|B|$.

An alternative approach to the Representation Theorems (i) and (ii)  \
for indefinite sesquilinear forms has been developed in
\cite{Fleige:1}, \cite{Fleige:2}, \cite{Fleige:3}
by Fleige, Hassi, and de Snoo in the framework of the Krein space theory.
Their results extended the list of  criteria equivalent to \eqref{ravny}. In particular,
it has been proven that the  condition \eqref{ravny} holds if and only if infinity is not a singular critical point (see \cite{Langer} for a discussion of this notion) for the range restriction $B_{\fa}$ in the Krein space $(\fK, [\cdot, \cdot ])$ of the operator $B$. Here
$\fK:=\Dom(A^{1/2})$, $[x,y]:=\fb[x,y]$ is an indefinite inner product on $\fK$, and
\begin{equation}\label{restriction}
 B_{\fa}=B|_{\fK} \quad \text{ on }\Dom  ( B_{\fa})=\{x\in \Dom (B)\subset \fK\, |\,  Bx\in \fK\}.
\end{equation}
A number of necessary and sufficient conditions for the regularity of the critical
point infinity have been discovered by \'Curgus in \cite{Curgus}.

The existence of operators in a Krein space with a singular
critical point at infinity is established in \cite{Fleige:0}, \cite{Fleige:00} and
\cite{Volkmer}. We remark that, by Proposition 5.3 in \cite{Fleige:2}, the existence of such operators implies the existence of a Hilbert space $\fH$ and a symmetric sesquilinear form $\fb$ on it such
that the condition \eqref{ravny} does not hold for the associated (by the First Representation Theorem) operator $B$ (see Example 5.4 in \cite{Fleige:2}). In this context, as a by-product of our considerations,
the range restriction $B_\fa$ of the operator $B$ in the Krein space $(\fK, \fb[\cdot, \cdot])$
constructed in Example \ref{exam} below provides a new fairly simple example
of an operator having infinity as a critical singular point (see Remark \ref{rem:m}).

%%%%%%%%%%%%%%%%%%%%%%%%%%%%%%%%%%%%%%%%%%%%%%%%%%%%%%%%%%%%%%%%%%%%%%%%%%%%%%%%%%%%%%%%%%%%%%%%%%%%%%%%%%%%%%%%%%%%%%%
\subsection*{Acknowledgments} The authors are grateful to anonymous referees for constructive comments on an earlier version of this paper. K.A.M.~is indebted to the Institute for Mathematics for
its kind hospitality during his two months stay at the Johannes Gutenberg-Universit\"{a}t Mainz in
the Summer of 2009. The work of K.A.M.~has been supported in part by the Deutsche Forschungsgemeinschaft,
Grant KO 2936/3-1, and by the Inneruniversit\"{a}ren For\-schungsf\"{o}rderung of the Johannes Gutenberg-Universit\"{a}t Mainz.
L.G.~has been supported by the exchange program between the University of Zagreb and the
Johannes Gutenberg-Universit\"{a}t Mainz and in part by the Grant
 037-0372783-2750 of the MZO\v{S}, Croatia. K.V.~has
been partly supported by the National Foundation of Science, Higher Education and Technical Development of the Republic of Croatia 2007-2009.
%%%%%%%%%%%%%%%%%%%%%%%%%%%%%%%%%%%%%%%%%%%%%%%%%%%%%%%%%%%%%%%%%%%%%%%%%%%%%%%%%%%%%%%%%%%%%%%%%%%%%%%%%%%%%%%%%%%%%%

%%%%%%%%%%%%%%%%%%%%%%%%%%%%%%%%%%%%%%%%%%%%%%%%%%%%
\section{Representation Theorems}\label{sec:repr}
%%%%%%%%%%%%%%%%%%%%%%%%%%%%%%%%%%%%%%%%%%%%%%%%%%%%

\begin{hypothesis}\label{hyp1}
Assume that $A$ and $H$ are self-adjoint operators in the Hilbert space $\fH$. Suppose that
\begin{itemize}
\item[(i)] $\inf \spec(A)>0$;
\item[(ii)] $H$ is bounded and has a bounded inverse;
\item[(iii)] the open interval $(h_-, h_+)$ is a maximal spectral gap of the operator $H$ containing $0$.
\end{itemize}
\end{hypothesis}

The following lemma introduces a  self-adjoint operator  naturally associated with the operators $A$ and $H$ from Hypothesis \ref{hyp1}.
\begin{lemma}\label{lem:1}
Assume Hypothesis \ref{hyp1}. Then the operator
\begin{equation}\label{operB}
B:=A^{1/2}HA^{1/2}
\end{equation}on the domain
$$\Dom(B)=\{x\in\Dom(A^{1/2})\,|\, HA^{1/2}x\in\Dom(A^{1/2})\}$$ is self-adjoint with a bounded inverse.
\end{lemma}

\begin{proof}
Introducing the bounded self-adjoint operator
\begin{equation*}
S := A^{-1/2} H^{-1} A^{-1/2},
\end{equation*}
we observe that $S$ has a trivial kernel and, hence, its inverse is a self-adjoint operator. It remains to note that $S^{-1}=B$.
\end{proof}

%%%%%%%%%%%%%%%%%%%%%%%%%%%%%%%%%%%%%%%%%%%%%%%%%%%%%%%%%%%%%%%%%%%%%%
\subsection{The First Representation Theorem}
%%%%%%%%%%%%%%%%%%%%%%%%%%%%%%%%%%%%%%%%%%%%%%%%%%%%%%%%%%%%%%%%%%%%%%

Under Hypothesis \ref{hyp1} consider the symmetric ses\-qui\-linear
form $\fb$ on $\Dom[\fb] = \Dom  (A^{1/2})$ defined by
\begin{equation}\label{form}
\fb[x,y]=\langle A^{1/2}x,HA^{1/2}y\rangle,\quad x,y\in \Dom[\fb] = \Dom  (A^{1/2}).
\end{equation}

\begin{theorem}\label{thm:main}
Assume Hypothesis \ref{hyp1} and  suppose that $\fb$ is
the symmetric sesquilinear form given by \eqref{form}.

Then the  operator $B$ referred to in Lemma \ref{lem:1} is the unique self-adjoint operator associated with the form $\fb$ in the sense that
\begin{equation}\label{eq:first}
\fb[x,y]= \langle x, By\rangle\quad\text{for all}\quad x\in\Dom[\fb],\quad y\in\Dom(B)\subset\Dom[\fb],
\end{equation}
with $\Dom(B)$ a core for $A^{1/2}$. Moreover,  the open interval
$(\alpha h_-, \alpha h_+)$, with $\alpha= \min \spec(A)$,  belongs to the  resolvent set of $B$.
\end{theorem}

\begin{proof}
By Lemma \ref{lem:1} the operator $B=A^{1/2}HA^{1/2}$
on $$\Dom(B)=\{x\in\Dom(A^{1/2})\,|\, HA^{1/2}x\in\Dom(A^{1/2})\}$$ is self-adjoint. It follows  that
\begin{equation*}
\fb[x,y] = \langle A^{1/2} x, H A^{1/2}y\rangle = \langle  x, A^{1/2} H A^{1/2}y\rangle
= \langle  x, By\rangle
\end{equation*}
for all $x\in\Dom(A^{1/2})$, $y\in\Dom(B)$, thereby proving the representation \eqref{eq:first}.

To prove that $\Dom(B)$ is a core for $A^{1/2}$ we assume
that $\langle y, A^{1/2}x\rangle_{\fH} = 0$ for some $y\in \fH$ and
for all $x\in\Dom(B)$. Since $\Dom(B)=\Ran(A^{-1/2}H^{-1}A^{-1/2})
$, one arrives at the conclusion that
$\langle y, H^{-1} A^{-1/2} z\rangle = 0$ for all $z\in\fH$.
Thus, $y=0$, since $\Ran(A^{-1/2})=\Dom(A^{1/2})$ is dense in $\fH$ and $H$ is an isomorphism.
Hence,  $\Dom(B)$
is a core for $A^{1/2}$.

Now we turn to the proof of the uniqueness. Assume that there exists a self-adjoint
operator $B^\prime$ with
 $\Dom(B^\prime)\subset\Dom[\fb]$ such that
\begin{equation*}
\langle x, B^\prime y\rangle = \fb[x,y]\quad\text{for all}\quad
x\in\Dom(\fb), \quad y\in\Dom(B^\prime).
\end{equation*}
Then
\begin{equation*}
\langle x, B^\prime y\rangle =
\fb[x,y]=\overline{\fb[y,x]}=\overline{\langle y, Bx\rangle} = \langle B x,
y\rangle
\end{equation*}
holds for all $x\in\Dom(B)$ and $y\in\Dom(B^\prime)$ which means that
$B^\prime$ is a restriction of $B^\ast$. Since both $B^\prime $
and $B$ are self-adjoint, we get $B^\prime = B$.

To complete the proof of the theorem it remains to show that the open interval
$(\alpha h_-,\alpha h_+)\ni 0$ belongs to the resolvent set of the operator $B$.
To this end we consider a  family of shifted quadratic forms
\begin{equation*}
\fb_\lambda [x,y]:=\fb[x,y]-\lambda \langle x,y\rangle=
\langle A^{1/2}x, (H-\lambda A^{-1})A^{1/2}x\rangle,\qquad
\lambda\in (\alpha h_-,\alpha h_+),
\end{equation*}
with $\Dom[\fb_\lambda]=\Dom[\fb]$. Observe that $H_\lambda := H-\lambda A^{-1}$ is bounded and has a bounded inverse if $\lambda\in(\alpha h_-,\alpha h_+)$. Indeed, the second resolvent
identity implies that
\begin{equation*}
H_\lambda^{-1} = H^{-1} + \lambda H^{-1} A^{-1/2} (I-\lambda A^{-1/2} H^{-1} A^{-1/2})^{-1}A^{-1/2}H^{-1}
\end{equation*}
holds as long as $I-\lambda A^{-1/2} H^{-1} A^{-1/2}$ is boundedly invertible. If the open interval $(h_-,h_+)\ni 0$ belongs to the resolvent set of the operator $H$, then
\begin{equation*}
h_-^{-1} I \leq H^{-1} \leq h_+^{-1} I.
\end{equation*}
Hence, we obtain the following bounds:
\begin{equation*}
\langle x, A^{-1/2} H^{-1} A^{-1/2}x\rangle = \langle A^{-1/2}x, H^{-1} A^{-1/2}x\rangle \leq h_+^{-1} \|A^{-1/2}x\|^2 \leq \frac{1}{\alpha h_+}
\end{equation*}
and
\begin{equation*}
-\langle x, A^{-1/2} H^{-1} A^{-1/2}x\rangle = -\langle A^{-1/2}x, H^{-1} A^{-1/2}x\rangle \leq - h_-^{-1} \|A^{-1/2}x\|^2 \leq -\frac{1}{\alpha h_-}.
\end{equation*}
Combining these bounds we arrive at the following two-sided  operator inequality
\begin{equation*}
\bigg (1+\frac{\lambda}{\alpha h_-}\bigg ) I\leq I-\lambda A^{-1/2} H^{-1} A^{-1/2} \leq \bigg (1-\frac{\lambda}{\alpha h_+}\bigg )I,
\end{equation*}
which shows that the operator $I-\lambda A^{-1/2} H^{-1} A^{-1/2}$ is boundedly invertible whenever $\lambda\in (\alpha h_-,\alpha h_+)$.

By the preceding arguments, there is a unique self-adjoint boundedly invertible
operator $B_\lambda$ with $\Dom(B_\lambda)\subset\Dom[\fb]$ such that
\begin{equation*}
\langle x, B_\lambda y\rangle = \fb_\lambda[x,y]\quad\text{for all}\quad x\in\Dom[\fb],\quad y\in\Dom(B_\lambda).
\end{equation*}
Clearly, $B_\lambda=B-\lambda I$, $\lambda\in (\alpha h_-,\alpha h_+)$, and, hence, the interval $(\alpha h_-,\alpha h_+)$ belongs to the resolvent set of the operator $B$, for $B_\lambda$ with $\lambda\in (\alpha h_-,\alpha h_+)$ has a bounded inverse. This completes the proof.
\end{proof}

\begin{remark}\label{rem:1}
Note that the operator $B$  referred to in Theorem  \ref{thm:main} is $A$-form bounded,  that is,
\begin{equation}\label{fcl}
\overline{A^{-1/2}BA^{-1/2}}=H\in\cB(\fH),
\end{equation}
where the bar denotes the closure of the operator
$A^{-1/2}BA^{-1/2}$ defined on
\begin{equation*}
\cD := A^{1/2}(\Dom(B)).
\end{equation*}
Indeed,  $A^{-1/2}BA^{-1/2}$ defined on $\cD$ coincides with $H|_\cD$.
Since $H$ is an isomorphism, from the representation
\begin{equation*}
A^{1/2} B^{-1}  = A^{1/2} A^{-1/2} H^{-1} A^{-1/2}  = H^{-1} A^{-1/2}
\end{equation*}
follows that the set $\cD$ is dense in $\fH$. Therefore, \eqref{fcl} holds.
\end{remark}

%%%%%%%%%%%%%%%%%%%%%%%%%%%%%%%%%%%%%%%%%%%%%%%%%%%%%%%%%%%%%%%%%%%%%%%%%%%%%%%%%%%%
\subsection{Applications to the case of off-diagonal form perturbations}
%%%%%%%%%%%%%%%%%%%%%%%%%%%%%%%%%%%%%%%%%%%%%%%%%%%%%%%%%%%%%%%%%%%%%%%%%%%%%%%%%%%%

\begin{theorem}\label{off}
Let $\fa$ be a positive definite closed symmetric sesquilinear form on
$\Dom[\fa]$ in a Hilbert space $\fH$  with the greatest lower bound $\alpha>0$ and
$\fv$ a symmetric $\fa$-bounded form on $\Dom[\fv]\supset\Dom[\fa]$, that is,
\begin{equation}\label{ogr}
v=\sup_{0\ne x\in \Dom[\fa]}\frac{|\fv[x]|}{\fa[x]}<\infty.
\end{equation}
Let $A$ be the uniformly positive self-adjoint operator associated with the closed form $\fa$ and
$J$ a self-adjoint involution commuting with $A$.

Assume, in addition, that the form $\fv$ is off-diagonal with respect to the orthogonal decomposition
\begin{equation*}
\fH=\fH_+\oplus\fH_-
 \quad \text{with}\quad
\fH_\pm=\Ran (I\pm J)\fH
\end{equation*}
in the sense that
\begin{equation*}
\fv[Jx,y]=-\fv[x,Jy], \quad x,y\in \Dom[\fa].
\end{equation*}
On $\Dom[\fb]:=\Dom[\fa]$  introduce the symmetric form
\begin{equation*}
\fb[x,y]=\fa[x,Jy]+\fv[x,y], \quad x,y\in \Dom[\fb].
\end{equation*}
Then there is a unique self-adjoint operator $B$ in $\fH$ such that $\Dom(B)\subset\Dom[\fb]$ and
\begin{equation*}
\fb[x,y]= \langle x, By\rangle\quad\text{for all}\quad x\in\Dom[\fb],\quad y\in \Dom (  B).
\end{equation*}
Moreover, the operator $B$ is boundedly invertible and the open interval
$(-\alpha, \alpha)\ni 0$ belongs to its resolvent set.
\end{theorem}

\begin{proof}
Due to the   hypothesis \eqref{ogr}, from the definition of the form $\fb$
follows that the sesquilinear form
\begin{equation}\label{H:def0}
\fh[x,y]:=\fb [A^{-\frac12}x, A^{-\frac12}y]
\end{equation}
with $\Dom[\fh]=\fH$ is bounded and symmetric.  Denote by $H$ the bounded self-adjoint operator associated with the form $\fh$.

Since the form $\fv$ is off-diagonal, the operator $H$ can be represented as the following $2\times 2$ block operator matrix
\begin{equation}\label{H:def}
H = \begin{pmatrix} I & T \\ T^\ast & -I \end{pmatrix},\quad T\in\cB(\fH_-,\fH_+),
\end{equation}
with respect to the orthogonal decomposition $\fH=\fH_+\oplus\fH_-$.

It is well know (see, e.g., \cite[Lemma 1.1]{KMM}
or \cite[Remark 2.8]{KMM:2})
that $H$ has a bounded inverse and, moreover, the open interval $(-1,1)$ belongs
 to the resolvent set of $H$. Thus, the operators $A$ and $H$ satisfy Hypothesis \ref{hyp1} and the claim follows by applying Theorem  \ref{thm:main}.
\end{proof}

\begin{remark}
Denote by $m_\pm$ the greatest lower bound of the form $\fa$ on the subspace $\fH_\pm$
(note that $\alpha=\min\{m_-,m_+\}$). If $m_-\ne m_+$ one can state that
not only the interval $(-\alpha, \alpha)$ but also the open
interval $(-m_-, m_+)\ni 0$ belongs to the resolvent set of the operator $B$.

Indeed, let
\begin{equation*}
\fa_\mu[x,y]-\mu\langle x, Jy\rangle, \quad x, y\in \Dom[\fa], \quad \mu\in (-m_-,m_+).
\end{equation*}
Then the form $\fa_\mu$ is closed and positive definite with
\begin{equation*}
v_\mu=\sup_{0\ne x\in \Dom[\fa]}\frac{|\fv[x]|}{\fa_\mu[x]}<\infty.
\end{equation*}
It is easy to see that the operator $B-\mu I$ is associated with the form
\begin{equation*}
\fb_\mu[x,y]=\fa_\mu[x,Jy]+\fv[x,y], \quad x,y\in \Dom[\fb].
\end{equation*}
Applying Theorem  \ref{off} one concludes that $B-\mu I$ has a bounded inverse for $\mu \in (-m_-, m_+)$ and hence
the open interval $(-m_-, m_+)\ni 0$ belongs to the resolvent set of the operator $B$.
\end{remark}

%%%%%%%%%%%%%%%%%%%%%%%%%%%%%%%%%%%%%%%%%%%%%%%%%%%%%%%%%%%%%
\subsection{A representation theorem for coercive forms}
%%%%%%%%%%%%%%%%%%%%%%%%%%%%%%%%%%%%%%%%%%%%%%%%%%%%%%%%%%%%%

Our next result shows how the self-adjoint operators $A$ and $H$ satisfying Hypothesis \ref{hyp1} naturally arise
in the context of perturbation theory.

\begin{theorem}\label{cor:1}
Assume that  $\fa$ is a positive definite closed symmetric ses\-qui\-linear form.
Let $A$ be the associated self-adjoint operator.
Suppose that $\fb$ is a symmetric $\fa$-bounded coercive sesquilinear form on $\Dom[\fb]=\Dom[\fa]$,
that is,
\begin{equation}\label{ner}
|\fb[x,y]| \leq \beta\, \fa[x,y]\quad\text{for all}\quad x,y\in \Dom[\fb],
\end{equation}
and
\begin{equation}\label{ner1}
|\fb[x,x]|\ge \alpha \fa[x,x] \quad \text{for all}\quad x\in \Dom[\fb],
\end{equation}
for some $0<\alpha\leq \beta$.

Then there is a unique bounded and boundedly invertible self-adjoint operator $H$ such that
the form $\fb$ admits the representation
\begin{equation}\label{ravH}
\fb[x,y]=\langle A^{1/2}x, H A^{1/2}y\rangle,\quad x,y\in \Dom[\fb] = \Dom  (A^{1/2}).
\end{equation}
\end{theorem}

\begin{proof}
{}From \eqref{ner}  follows that there exists a
bounded self-adjoint operator $M$ on the
Hilbert space $\fH_{\fa}:=\Dom[\fa]$ equipped with the inner product $\fa[\cdot,\cdot]$
such that \begin{equation}\label{ravM}
\fb[x,y]=\fa[x,My]=\langle A^{1/2}x, A^{1/2} My\rangle, \quad x, y\in \Dom(A^{1/2}).
\end{equation}
On the other hand, from  \eqref{ner1} one concludes that $M$ has a bounded inverse
with the bound
$\|M^{-1}\|\leq 1/\alpha$
(this is a special case of the classical Lax-Milgram lemma, see, e.g., \cite[Theorem IV.1.2]{Evans}).

By Lemma VI.3.1 in \cite{Kato} there exists a bounded self-adjoint operator $H$ in $\fH$ such that
\begin{equation}\label{eq:eq:1}
\fb[x,y]=\langle A^{1/2}x, H A^{1/2}y\rangle, \quad x, y\in \Dom(A^{1/2}).
\end{equation}

Comparing \eqref{ravM} and \eqref{eq:eq:1} yields the equality
\begin{equation*}
H A^{1/2}y=A^{1/2}My\quad\text{for all}\quad y\in\Dom(A^{1/2})
\end{equation*}
and, therefore,
\begin{equation*}
H=A^{1/2} M A^{-1/2}\quad \text{on}\quad \fH.
\end{equation*}

Since $M$ has a bounded inverse in $\fH_{\fa}$, it is an isomorphism of $\fH_{\fa}$.
It remains to note that $A^{-1/2}$ maps $\fH$ onto $\fH_{\fa}$ isomorphically and, therefore, $H$ is an isomorphism of $\fH$ and hence the self-adjoint operator $H$ has a bounded inverse.

The proof is complete.
\end{proof}

\begin{remark}\label{rem:mcintosh}
According to McIntosh \cite{McIntosh:1}, a possibly
sign-indefinite sesquilinear form that admits the representation
$\fb[x,y]=\fa[x, My]$, with $M$ an isomorphism of $\fH_\fa$
(cf., e.g., \eqref{ravM}), is called $0$-closed.
We refer to Theorem 3.2 in \cite{McIntosh:1} to emphasize the role of
$\,0$-closed forms that they play in the context of the First Representation Theorem.
\end{remark}

\begin{remark}\label{rem:fleige}
An alternative notion of closedness of indefinite quadratic forms in the Krein space setting has been introduced in
\cite{Fleige:1}, \cite{Fleige:2}, \cite{Fleige:3}
by Fleige, Hassi, and de Snoo. The sesquilinear form
$\fb$ referred to in Theorem \ref{cor:1} is closed in this sense as well. In other words,
$(\Dom(A^{1/2}),$ $\fb[\cdot,\cdot])$ is a Krein space continuously embedded in $\fH$.

Indeed, since the inner products on $\fH_{\fa}$ given by
$$\langle x,y\rangle_{\fa} = \langle A^{1/2}x, A^{1/2}y\rangle$$
and by $$\langle x,y\rangle_H=\langle A^{1/2}x, |H| A^{1/2}y\rangle$$ are equivalent,
the Hilbert space $\fH_{\fa}$ equipped with the new inner product $\langle \cdot,\cdot\rangle_H$ is continuously embedded in $\fH$. The involution $J:=A^{-1/2} \sign(H) A^{1/2}$ is bounded in $\fH_{\fa}$, for
\begin{equation*}
 |\langle x, Jy \rangle_{\fa}| = |\langle A^{1/2}x, A^{1/2} A^{-1/2}\sign(H)A^{1/2}x \rangle|\leq \|x\|_{\fa} \|y\|_{\fa}
\end{equation*}
holds for all $x,y\in\fH_{\fa}$. Since
\begin{equation*}
 \langle x, Jx\rangle_H = \langle A^{1/2}x, |H|A^{1/2}x\rangle
\end{equation*}
is real, $J$ is self-adjoint in $\fH_{\fa}$ with respect to the inner product
$\langle \cdot, \cdot\rangle_H$. Taking into account the equality
\begin{equation*}
 \fb[x,y] = \langle A^{1/2}x, H A^{1/2}y\rangle = \langle x, Jy\rangle_H,\quad x,y\in\fH_{\fa},
\end{equation*}
one concludes that $(\Dom(A^{1/2}),\fb[\cdot,\cdot])$ is a Krein space.
\end{remark}

%%%%%%%%%%%%%%%%%%%%%%%%%%%%%%%%%%%%%%%%%%%%%%%%%%%%%%%%%%%%%%%
\subsection{The Second Representation Theorem}
%%%%%%%%%%%%%%%%%%%%%%%%%%%%%%%%%%%%%%%%%%%%%%%%%%%%%%%%%%%%%%%

Our next goal is to prove the Second Representation Theorem (ii)
 for sign-indefinite forms, a result
 originally due to McIntosh \cite{McIntosh:1}.
We emphasize that equation  \eqref{eq:2RT} combined with
\eqref{eq:1RT}   shows that not only the operator $B$ is associated with the form $\fb$ (Theorem \ref{thm:main}), but also that the form $\fb$ is represented by the operator $B$, provided
that the form domain stability condition
 \begin{equation}\label{formdom}
\Dom(A^{1/2})=\Dom(|B|^{1/2})
\end{equation}
holds.

More precisely, we have the following result.

\begin{theorem}\label{thm:main:2}
Assume hypotheses of Theorem \ref{thm:main} and let $B$ be the operator referred to therein.

If $\Dom(A^{1/2})=\Dom(|B|^{1/2})$, then
\begin{equation}\label{eq:second}
\fb[x,y]= \langle |B|^{1/2} x, \sign(B)|B|^{1/2} y\rangle\quad\text{for all}\quad x,y\in\Dom[\fb]=\Dom(|B|^{1/2}).
\end{equation}
\end{theorem}

\begin{proof}
{}From Theorem \ref{thm:main} follows that
\begin{equation*}
\fb[x,y]=\langle  x, B y\rangle
\end{equation*}
for all $x\in\Dom(|B|^{1/2})$, $y\in\Dom(B)$, which yields
\begin{equation}\label{eq:prep}
\fb[x,y]=\langle  |B|^{1/2} J x, |B|^{1/2}  y\rangle
\quad \text{for all } \quad x\in\Dom(|B|^{1/2}),\,\,\, y\in\Dom(B),
\end{equation}
where $J=\sign B$.

To complete the proof it remains to show that \eqref{eq:prep}
 holds for all $x,y\in\Dom(|B|^{1/2})$. To this
end we fix $x\in\Dom(|B|^{1/2})$ and consider two linear functionals given by
\begin{equation*}
\begin{split}
\ell_1(y) &:= \fb[x,y],\\
\ell_2(y) &:= \langle  |B|^{1/2} J x, |B|^{1/2}  y\rangle
\end{split}
\end{equation*}
defined on $\Dom(A^{1/2})\equiv \Dom(|B|^{1/2})$. For the form $\fb$ is $\fa$-bounded, the functional $\ell_1$ is continuous on $\fH_{\fa}$. Since $\Dom(A^{1/2})=\Dom(|B|^{1/2})$, by the closed graph theorem the operator
$A^{1/2}$ is $|B|^{1/2}$-bounded and $|B|^{1/2}$ is $A^{1/2}$-bounded. Therefore the norms
on $\fH_{\fa}$
\begin{equation*}
|x|_{\fa} := \|A^{1/2}x\|\quad\text{and}\quad |x|_{\fb} := \||B|^{1/2}x\|
\end{equation*}
are equivalent. The functional $\ell_2$ is continuous on $\Dom(|B|^{1/2})$ in the topology of the norm $|\cdot|_{\fb}$ and, thus, it is continuous on $\fH_{\fa}$.
 Since $\Dom(B)$ is a core for $|B|^{1/2}$, it follows that $\Dom(B)$ is dense in $\fH_{\fa}$. Hence, since by \eqref{eq:prep} the functionals $\ell_1$ and $\ell_2$ agree on the set $\Dom(B)$ dense in $\fH_{\fa}$, it follows that $\ell_1=\ell_2$ on $\fH_{\fa}$.
\end{proof}

Under the form domain stability condition \eqref{formdom}, Theorem \ref{thm:main:2}
combined with Theorem \ref{thm:main} establishes a one-to-one correspondence between
the symmetric forms of the type \eqref{form} and the associated self-adjoint operators $B$ given by \eqref{operB}.

The following example provides a pair of self-adjoint
operators $A$ and $H$  satisfying  Hypothesis
\ref{hyp1} such that the form domain stability condition
\eqref{formdom} required in the hypothesis of Theorem \ref{thm:main:2} does not hold.

\begin{example}\label{exam}
In the Hilbert space $\displaystyle\fH = \ell^2(\N;\C^2) \cong \ell^2(\N)\oplus\ell^2(\N)$ consider the self-adjoint operator
\begin{equation}\label{A:def}
A=\bigoplus_{k\in\N} \begin{pmatrix} 1 & 0 \\ 0 & k^2 \end{pmatrix}\quad \text{on}\quad  \Dom(A)=\ell^{2,0}(\N)\oplus \ell^{2,4}(\N),
\end{equation}
where $\ell^{2,p}(\N)$ denotes the space of sequences $\{a_k\}_{k=1}^\infty$ such that $\sum_{k\in \N}k^p|a_k|^2<\infty$,
and the bounded self-adjoint operator $H$ given by
\begin{equation*}
H=\bigoplus_{k\in\N} \begin{pmatrix} 0 & 1 \\ 1 & 0 \end{pmatrix} .
\end{equation*}

A simple computation shows that the operator $B=A^{1/2}HA^{1/2}$ associated with the form \eqref{form}
admits the representation
\begin{equation*}
 B  = \bigoplus_{k\in\N} \begin{pmatrix} 0 & k \\ k & 0 \end{pmatrix}
\quad \text{ on } \quad\Dom(B)=\ell^{2,2}(\N)\oplus\ell^{2,2}(\N).
\end{equation*}
Hence,
\begin{equation*}
|B|= \bigoplus_{k\in\N} \begin{pmatrix} k & 0 \\ 0 & k \end{pmatrix}\quad\text{ on }\quad \Dom(|B|)=\Dom(B).
\end{equation*}
Clearly,
\begin{equation*}
\Dom(A^{1/2})=\ell^{2, 0}(\N)\oplus \ell^{2, 2}(\N)\quad \text{ and }\quad
\Dom(|B|^{1/2})=\ell^{2, 1}(\N)\oplus \ell^{2, 1}(\N)
\end{equation*}
and, therefore,
\begin{equation}\label{ungl}
\Dom(|B|^{1/2})\neq \Dom (A^{1/2}).
\end{equation}
\end{example}

In the particular case considered in Example \ref{exam} we face the following
phenomenon, which apparently never happens whenever $B$ is semi-bounded (cf.~Lemma \ref{thm:main:20}
below): The self-adjoint operator $B$ is associated with two different sesquilinear forms
$\fb_1$ and $\fb_2$ given by
\begin{align*}
\fb_1[x,y]&=\langle A^{1/2}x,HA^{1/2}y\rangle,\quad
x,y\in \Dom[\fb] = \Dom  (A^{1/2}),\\
\fb_2[x,y]&= \langle |B|^{1/2} x, \sign(B)|B|^{1/2} y\rangle,\quad
x,y\in\Dom[\fb]=\Dom(|B|^{1/2}),
\end{align*}
but only one of them, namely the form $\fb_2$, is represented by $B$. We will turn back to the discussion of this
phenomenon in Section \ref{sec:conv}.

\begin{remark}\label{rem:m}
Let $B_{\fa}$ denote the range restriction of the operator
$B$ from Example \ref{exam} in the Krein space
$(\fK,[\cdot,\cdot])$,  with $\fK=\Dom(A^{1/2})$,
the indefinite inner product
\begin{equation*}
[x,y]=\langle A^{1/2} x, H A^{1/2}y\rangle_{\fH},
\end{equation*}
and the fundamental symmetry
\begin{equation*}
J := A^{-1/2} H A^{1/2} = \bigoplus_{k\in\N}\begin{pmatrix} 0 & k \\ k^{-1} & 0 \end{pmatrix}
\end{equation*}
(cf.~\eqref{restriction} and Remark \ref{rem:fleige}). One easily verifies that
\begin{equation*}
\Dom(B_{\fa}) = \ell^{2,4}(\N) \oplus \ell^{2,2}(\N)
\end{equation*}
and that
\begin{equation*}
\begin{split}
[x, B_\fa x] &= \langle A^{1/2} x, H A^{1/2} B_{\fa} x \rangle_{\fH}\\
&= \langle A^{1/2} H A^{1/2} x, A^{1/2} H A^{1/2} x \rangle_{\fH} = \|Bx\|_{\fH}^2
\end{split}
\end{equation*}
for all $x\in\Dom(B_{\fa})$, that is, the operator $B_{\fa}$ is positive with respect to the indefinite inner product $[\cdot,\cdot]$. Since $B$ is boundedly invertible, its range restriction $B_{\fa}$ is boundedly invertible as well. Hence, according to Theorem 2.5 in \cite{Curgus}, infinity is not a singular critical point of the operator $B_{\fa}$ if and only if the norms generated by the positive definite inner products $[\cdot,J\cdot]$ and $[\cdot, B_{\fa}|B_{\fa}^{-1}|\cdot]$ on $\Dom(JB_{\fa})=\Dom(B_{\fa})$ are equivalent. Since
\begin{equation*}
[x, Jy] = \langle A^{1/2}x, A^{1/2}y\rangle_{\fH}\quad \text{and} \quad
[x, B_{\fa}|B_{\fa}^{-1}|y]  = \langle |B|^{1/2}x, |B|^{1/2}y\rangle_{\fH}
\end{equation*}
hold for all $x,y\in\Dom(B_{\fa})$, the corresponding  norms are equivalent if and only if the domains $\Dom(A^{1/2})$ and $\Dom(|B|^{1/2})$ agree.
Therefore, due to \eqref{ungl}, infinity is a singular critical point of the operator $B_{\fa}$.
\end{remark}

%%%%%%%%%%%%%%%%%%%%%%%%%%%%%%%%%%%%%%%%%%%%%%%%%%%%%%%%%%%%%%%
\section{Form-domain Stability Criteria}\label{sec:stability}
%%%%%%%%%%%%%%%%%%%%%%%%%%%%%%%%%%%%%%%%%%%%%%%%%%%%%%%%%%%%%%%

The main goal of this section is to establish a number of criteria ensuring the form-domain
stability condition \eqref{formdom}.

The following simple functional-analytic lemma plays a key role in our further considerations.

\begin{lemma}\label{lem:3}
Let $(\fH,\langle\cdot,\cdot\rangle)$ and $(\fH',\langle\cdot,\cdot\rangle')$ be Hilbert spaces. Assume that $\fH'$ is continuously embedded in $\fH$. If $T$ is a continuous linear map from $\fH$ to $\fH$ leaving $\fH'$ invariant (as a set), then the operator $T'$ induced by $T$ on $\fH'$ is continuous (in the topology of $\fH'$).
\end{lemma}

\begin{proof}
By the hypothesis of the lemma the operator  $T'$ is defined on the whole of $\fH'$. Therefore, by the closed graph theorem it suffices to prove that $T'$ is closed. Assume that
\begin{equation*}
x_n \overset{\fH'}{\longrightarrow}x\qquad\text{and}\qquad T' x_n \overset{\fH'}{\longrightarrow}y.
\end{equation*}
Since the Hilbert space $\fH'$ is continuously embedded into $\fH$, one also has
\begin{equation*}
x_n \overset{\fH}{\longrightarrow}x\qquad\text{and}\qquad T x_n \overset{\fH}{\longrightarrow}y.
\end{equation*}
{}From the continuity of $T$ in $\fH$, it follows that $Tx=y$ in $\fH$, and, thus, $T' x=y$ in $\fH'$, which proves the claim.
\end{proof}

Introduce the following symmetric  nonnegative operators
\begin{equation}\label{op:W}
X:=A^{-1/2}|B|A^{-1/2}\quad \text{ on }\quad  \Dom(X)=A^{1/2}\Dom (B)
\end{equation}
and
\begin{equation}\label{op:M}
Y:=A^{1/2}|B|^{-1}A^{1/2}\quad \text{ on } \quad \Dom(Y)=\Dom(A^{1/2}).
\end{equation}

By Remark \ref{rem:1},  $\Dom(X)=A^{1/2} \Dom(B)$ is dense in $\fH$.
Hence, $X$ is a densely defined operator, so is $Y$, since $\Dom(A^{1/2})$ is obviously a dense set.

Now we are prepared to present the main result of this  section.

\begin{theorem}\label{lem:n} Let the operators $A$ and $B$ be as in Theorem \ref{thm:main}. Then
the following are equivalent:
\begin{itemize}
\item[(i)] $\Dom(|B|^{1/2})=\Dom(A^{1/2})$;
\item[(ii)] $\Dom(A^{1/2})\subset\Dom(|B|^{1/2})$;
\item[(ii$'$)] $\Dom(|B|^{1/2})\subset\Dom(A^{1/2})$;
\item[(iii)] $X=A^{-1/2}|B|A^{-1/2}$ is a bounded symmetric \
operator on $\Dom(X)=A^{1/2}\Dom (B)$;
\item[(iii$'$)] $Y=A^{1/2}|B|^{-1}A^{1/2}$ is a bounded symmetric operator on $\Dom(Y)=\Dom(A^{1/2})$,
\item[(iv)] $K:=A^{1/2} \sign(B)A^{-1/2}$ is a bounded involution on $\fH$;
\item[(v)] $\sign(B) \Dom(A^{1/2}) \subset \Dom(A^{1/2})$.
\end{itemize}
\end{theorem}

\begin{proof}
The implications (i)$\Rightarrow$(ii)  and  (i)$\Rightarrow$(ii$'$) are  obvious.

(ii)$\Rightarrow$(iii). Since $\Dom(A^{1/2})\subset\Dom(|B|^{1/2})$, the operator $|B|^{1/2}A^{-1/2}$  is bounded.
Introducing the sesquilinear form
\begin{equation}
\fx [x,y]=\langle
x, Xy \rangle \quad \text{ on } \quad \Dom[\fx]=\Dom (X),
\end{equation}
one concludes that the form $\fx$ can also be  represented as a bounded form (since
$|B|^{1/2}A^{-1/2}$   is bounded)
\begin{equation*}
\fx [x,y]=\langle
|B|^{1/2}A^{-1/2}x,|B|^{1/2}A^{-1/2} y \rangle.
\end{equation*}
Therefore, the sesquilinear form $\fx$ is associated with a
bounded operator and, hence,
 the closure of $X$ is a bounded operator
defined on the whole Hilbert space $\fH$.

(ii$'$)$\Rightarrow$(iii$'$). Arguing as above,
one shows that the operator $A^{1/2}|B|^{-1/2}$  is bounded and therefore the form
\begin{equation*}
\fy [x,y]=\langle
x, Yy \rangle= \langle
A^{1/2}|B|^{-1/2}x,A^{1/2}|B|^{-1/2} y \rangle\quad \text{ on } \quad \Dom[\fy]=\Dom (Y),
\end{equation*}
is a bounded form. Hence,  the closure of $Y$
is a bounded operator defined on the whole Hilbert space $\fH$.

(iii)$\Rightarrow$(iv).
%It is straightforward to verify that the restriction of the bounded operator $\left(|B|
 %^{1/2}A^{-1/2}\right)^\ast |B|^{1/2}A^{-1/2}$ onto $\cD:=A^{1/2}\Dom(B)\subset\fH$ is given by %$A^{-1/2} |B|A^{-1/2}$.
Note that the operator $K$ on its natural domain
\begin{equation*}
\Dom(K)=\{x\in\fH\, |\, \sign(B) A^{-1/2}x\in\Dom(A^{1/2})\}
\end{equation*}
is obviously closed.  Moreover, it is also clear that
\begin{equation*}
\Dom (X)=A^{1/2}\Dom(B)\subset \Dom(K).
\end{equation*}
Since for any $x\in \Dom(X)$ one gets that
\begin{align*}
H^{-1} X x&=H^{-1}A^{-1/2} |B| A^{-1/2}x=A^{1/2}A^{-1/2}H^{-1}A^{-1/2} |B| A^{-1/2}x\\
&=
A^{1/2}B^{-1}|B|A^{-1/2}x=A^{1/2}\sign (B)A^{-1/2}x=Kx
\end{align*}
and both $H^{-1}$ and $X$ are bounded operators, one concludes that $K|_{\Dom(X)}$
is a bounded operator.  Since $K$ is closed and $\Dom(X)\subset \Dom (K)$, the operator
$K|_{\Dom(X)}$ is closable. Since $K|_{\Dom(X)}$ is bounded, the domain of its closure is a closed subspace that contains
a set dense in $\fH$.
Therefore, $K=\overline{K}=\overline{K|_{\Dom(X)}}$ is a bounded involution with $\Dom(K)=\fH$.

(iii$'$)$\Rightarrow$(iv). For any $x\in H^{-1}\Dom(Y)=H^{-1}\Dom(A^{1/2})$ one gets that
\begin{align*}
YHx&=A^{1/2}|B|^{-1}A^{1/2}Hx=A^{1/2}|B|^{-1}A^{1/2}HA^{1/2}A^{-1/2}x\\
&=A^{1/2}|B|^{-1}BA^{-1/2}x=A^{1/2}\sign (B)A^{-1/2}x=Kx.
\end{align*}
Next we check that the dense set $H^{-1}\Dom(A^{1/2})$ is a subset of $\Dom (K)$.
Indeed, if $x\in H^{-1}\Dom(A^{1/2})$, then $x=H^{-1}A^{-1/2}y$ for some $y\in \fH$. Hence,
\begin{equation*}
\sign (B)A^{-1/2}x=\sign (B) A^{-1/2}H^{-1}A^{-1/2}y=\sign (B)B^{-1}y\in\Dom(B)\subset \Dom(A^{1/2}),
\end{equation*}
and, therefore, $x\in \Dom(K)$. Now, to conclude that $K$ is a bounded involution it remains to argue as in the proof of the implication (iii)$\Rightarrow$(iv).

(iv)$\Rightarrow$(v). {}Since the operator $K=A^{1/2} \sign(B)A^{-1/2}$ is well defined as a bounded operator
on the whole Hilbert space  $\fH$,  one cocnludes  that $\sign(B)$ must leave $\Dom(A^{1/2})$ invariant.

(v)$\Rightarrow$(i). We start with the particular case of positive $H$. Consider the positive definite sesquilinear form
\begin{equation*}
\fb[x,y] := \langle A^{1/2} x, H A^{1/2} y\rangle
\end{equation*}
defined on $\Dom[\fb]=\Dom(A^{1/2})$. Since $H$ is positive, one can represent the form $\fb$ as
\begin{equation*}
\fb[x,y] = \langle H^{1/2}A^{1/2} x, H^{1/2} A^{1/2} y\rangle.
\end{equation*}
For $H^{1/2} A^{1/2}$ is closed, the form $\fb$ is closed.

By the definition of the operator $B$ (cf.~Lemma \ref{lem:1}),
\begin{equation*}
\fb[x,y] = \langle  x, A^{1/2} H A^{1/2} y\rangle= \langle  x, B y\rangle
\end{equation*}
for all $x\in\Dom[\fb]\equiv\Dom(A^{1/2})$ and $y\in\Dom(B)$. By Lemma \ref{lem:1} the operator $B$ is self-adjoint. Therefore, the operator $B$ is associated with the form $\fb$. The second representation theorem for positive definite sesquilinear forms \cite{Kato} yields $\Dom[\fb]=\Dom(B^{1/2})$, which proves the claim.

We turn to the case when $H$ is not necessarily positive.

Set for brevity $J:=\sign(B)$. Denote by $J_\fa$ the operator on $\fH_\fa$ induced by $J$. Since $J^2=I$, by Lemma \ref{lem:3} the operator $J_\fa$ is a bounded involution, not necessarily unitary.

This observation allows one to conclude that
\begin{equation*}
K = A^{1/2} J A^{-1/2}
\end{equation*}
is a bounded involution in the Hilbert space $\fH$. To complete the proof of the theorem one
notices that
\begin{equation}\label{eq:modulb}
|B|=B J = A^{1/2} HK A^{1/2}.
\end{equation}
Since $|B|\geq 0$, one immediately verifies that the bounded operator $HK$ is nonnegative. Since both $H$ and $K$ are Hilbert space isomorphisms, the self-adjoint operator $HK$ has a bounded inverse. Since the case of positively definite $H$ has been already discussed, we arrive at the conclusion that $\Dom(|B|^{1/2})=\Dom(A^{1/2})$.
\end{proof}

\begin{remark}
We note that the equivalence of (i) and (iv) has been established by McIntosh in \cite[Lemma 2.5]{McIntosh:2}.
Independently, when our work has already appeared as a preprint
arXiv:1003.1908, the equivalence of
(i),(ii), (ii$'$) has been established by a different method in
\cite[Proposition 2.5]{Fleige:4}.

\end{remark}

\begin{remark}\label{rem:2}
If the domains $\Dom(|B|^{1/2})$ and $\Dom(A^{1/2})$ agree, then the sesquilinear form $|\fb|$ associated with the positive operator $|B|$ can be represented as
\begin{equation*}
|\fb|[x,y]= \langle (HK)^{1/2} A^{1/2}x, (HK)^{1/2} A^{1/2}y\rangle,\quad x,y\in\Dom(|B|^{1/2})=\Dom(A^{1/2})
\end{equation*}
and, therefore, along with \eqref{eq:modulb} one obtains the factorization
\begin{equation*}
|B| = \left( (HK)^{1/2} A^{1/2}\right)^\ast (HK)^{1/2} A^{1/2}.
\end{equation*}
\end{remark}

\begin{remark}\label{rem:4}
If the operator $K$ is bounded, then it is similar
to a unitary self-adjoint operator.
Indeed, the equalities $K^2=I$ and $K^{-1}=K$ imply that both
$K$ and $K^{-1}$ are power bounded, that is, $\sup_{n\in\Z}\|K^n\|$
is finite. Hence, by \cite[Theorem 2.4]{Casteren}, $K$ is
similar to a unitary operator $U$. Since $K^{-1}=K$,
 we have $U^{-1}=U=U^*$ such that $U$ is self-adjoint.
\end{remark}

%%%%%%%%%%%%%%%%%%%%%%%%%%%%%%%%%%%%%%%%%%%%%%%%%%
\subsection{Sufficient criteria}
%%%%%%%%%%%%%%%%%%%%%%%%%%%%%%%%%%%%%%%%%%%%%%%%%%

The following lemma provides several sufficient (but not necessary) criteria
for the form-domain stability condition to hold.

\begin{lemma}\label{thm:main:20}
Assume hypotheses of Theorem \ref{thm:main} and let $B$ be the operator referred to therein.
If one of the following conditions
\begin{itemize}
\item[(i)] the operator
$H$ maps $\Dom(A^{1/2})$ onto itself;
\item[(ii)] the operator $H$ is uniformly positive;
\item[(iii)] the operator $B$ is semi-bounded;
\end{itemize}
hold, then
\begin{equation}\label{formdom1}
\Dom (|B|^{1/2})=\Dom(A^{1/2}).
\end{equation}

\end{lemma}

\begin{proof}
(i). Since $\Dom(A^{1/2})$ is $H$-invariant, and
$$\Dom(B)=\{x\in\Dom(A^{1/2})\,|\, HA^{1/2}x\in\Dom(A^{1/2})\}$$
one concludes  that
 $\Dom(B)=\Dom(A)$,  which implies \eqref{formdom1} by the Heinz inequality (cf.~\cite[Theorem\ 3]{He51},
\cite[Theorem\ IV.1.11]{KPS82}, \cite[Ch. 10, Section 4]{Birman}).

(ii). If $H$ is uniformly positive, then the operator $B$ is nonnegative and, therefore, $B$ is also uniformly positive by the First Representation Theorem. Thus, $\sign(B)=I$ and condition (iv) of Theorem \ref{lem:n}
is trivially fulfilled and, hence, \eqref{formdom1} holds.

(iii). Assume, for definiteness, that the operator $B$ is semi-bounded from below,
and, hence, $B+\beta I\ge 0$ for  $\beta >|\inf\spec(B)|$. Therefore, for those $\beta$, one gets that
\begin{equation*}
B+\beta I=A^{1/2}HA^{1/2}+\beta I=A^{1/2}(H+\beta A^{-1})A^{1/2}
\end{equation*}
and, moreover, $H+\beta A^{-1}\ge 0$.

Since $H$, by hypothesis,  has a bounded inverse, by the Birman-Schwinger principle $H+\beta A^{-1}$ has a bounded inverse if and only if the operator $I+\beta A^{-1/2}H^{-1}A^{-1/2}=I+\beta B^{-1}$ does,
which is the case for $\beta >|\inf \spec(B)|$. Thus, $H+\beta A^{-1}$ is uniformly positive and by (ii) one obtains that
\begin{equation*}
\Dom |B+\beta I|^{1/2}=\Dom(A^{1/2}).
\end{equation*}
It remains to remark that
\begin{equation*}
\Dom |B+\beta I|^{1/2}=\Dom (B+\beta I)^{1/2}=\Dom(|B|^{1/2})
\end{equation*}
and the claim follows.
\end{proof}

%%%%%%%%%%%%%%%%%%%%%%%%%%%%%%%%%%%%%%%%%%%%%%%%%%%%%%%%%%%%%%%%%
\subsection{The form domain stability in pictures}
%%%%%%%%%%%%%%%%%%%%%%%%%%%%%%%%%%%%%%%%%%%%%%%%%%%%%%%%%%%%%%%%%

Given a not necessarily semibounded self-adjoint operator $A$, introduce the Sobolev-like
scale of spaces
\begin{equation*}
\cH_A^s=\Dom(|A|^{s/2}),\quad s\ge 0,
\end{equation*}
equipped with the graph norm of $|A|^{s/2}$, with a natural convention that $\cH_A^0=\fH$, the underlying Hilbert space.

We remark that if self-adjoint operators $A$ and $B$ have  bounded inverses,
and $\cH_A^s=\cH_B^s$ for some $s>0$, then$\cH_A^t=\cH_B^t$ for
all $0\le t\le s$ by the Heinz inequality.
In particular, under hypothesis (i) of Lemma \ref{thm:main:20},
the domains $\Dom(A)$ and $\Dom(B)$ of $A$ and $B$ coincide. That is,
$\cH_A^2=\cH_B^2$, and, therefore, the form domain stability
condition $\cH_A^1=\cH_B^1$ holds automatically.

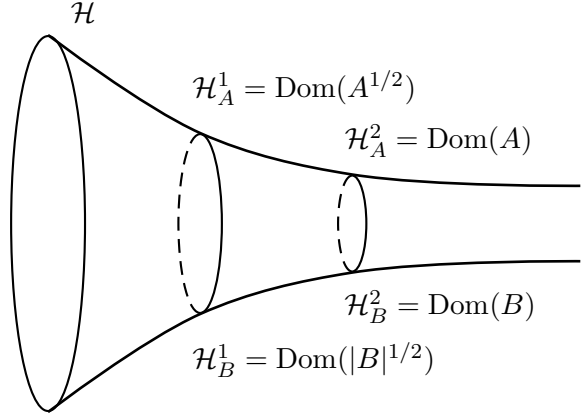
\begin{figure}[htb]
\centering
\begin{pspicture}[](0,-3)(7,3)
\pscurve[linecolor=black,linewidth=1pt]
(0,2.5)
(2,1.2)
(4,0.65)
(7,.5)
%\psecurve[linestyle=dashed](3,4)(4,4.23)(5,4.43)
\pscurve[linewidth=1pt]
(0,-2.5)
(2,-1.2)
(4,-0.65)
(7,-.5)

\psellipse[](0,0)(0.5,2.5)
\psellipticarc[linestyle=dashed](2,0)(0.3,1.2){90}{270}
\psellipticarc[](2,0)(0.3,1.2){-90}{90}

\psellipticarc[](4,0)(0.2,.65){-90}{90}
\psellipticarc[linestyle=dashed](4,0)(0.2,.65){90}{270}

%%%%%%%%%%%%%%%%%%%%%%%%%%%%%%%%%%%%%%%%%%%%%%%%%%%%%
%\psellipse[linecolor=red](1.58,-.36)(0.69,.6)
%\psellipticarc[linecolor=red](1.58,-.36)(1,.6){-100}{0}
%%%%%%%%%%%%%%%%%%%%%%%%%%%%%%%%%%%%%%%%%%%%%%%%%%%%%%%%
%%%%%%%%%%%%%%%%%%%%%%%%%%%%%%%%%%%%%%%%%%%%%%%%%%%%%%%%%
\uput{0.3}[0](0,2.8){$\cH$}
\uput{0.3}[0](1.6,1.8){$\cH^1_A=\Dom(A^{1/2})$}
\uput{0.3}[0](3.6,1.1){$\cH^2_A=\Dom(A)$}
\uput{0.3}[0](1.6,-1.8){$\cH^1_B=\Dom(|B|^{1/2})$}
\uput{0.3}[0](3.6,-1.1){$\cH^2_B=\Dom(B)$}
\end{pspicture}
\caption{\label{fig:1} The Sobolev-like scale of
spaces for the operators $A$ and $B$. $\cH_A^2=\cH_B^2$ and, hence,
the form domain stability condition
$\cH_A^1=\cH_B^1$ holds.}
\end{figure}

The diagram depicted in Fig.~\ref{fig:1} illustrates the case.

%\vskip .5cm

Under Hypothesis \ref{hyp1} the perturbation may change the domain of $A$, so that
$\Dom(A)\ne \Dom(B)$. However, the form domain stability condition may still hold.
For instance, it is the case when the operator  $B$ is semibounded.
A typical diagram is presented in Fig.~\ref{fig:2}.

\begin{figure}[htb]
\centering
%\vspace{20mm}
\begin{pspicture}[](0,-3)(7,3)
\pscurve[linecolor=black,linewidth=1pt]
(0,2.5)
(2,1.2)
(4,0.65)
(7,.5)
%\psecurve[linestyle=dashed](3,4)(4,4.23)(5,4.43)
\pscurve[linecolor=black,linewidth=1pt]
(0,-2.5)
(2,-1.2)
(4,-0.65)
(7,-.5)
%\psecurve[linestyle=dashed](3,4)(4,4.23)(5,4.43)
\psellipse[](0,0)(0.5,2.5)
\psellipticarc[linestyle=dashed](2,0)(0.3,1.2){90}{270}
\psellipticarc[](2,0)(0.3,1.2){-90}{90}
%\pscurve[linestyle=dashed,dash=.5pt 1pt](1.46,-1.48)(1.29,-1.44)(1.18,-1.28)(1.73,-0.50)(1.72,-.25)
%\pscurve[](1.48,-1.48)(1.7,-1.3)(2,-.95)(2.28,0.0)(2.27,0.50)(2.33,0.8)(2.43,.934)(2.65,0.93)
%\pscurve[linestyle=dashed,dash=.5pt 1pt](2.65,0.93)
%(2.7,.9)(2.72,.75)(2.71,.7)(2.7,0.65)(2.6,.55)(2.54, 0.52)(2.45, 0.51)(2.3, 0.5)(2.1, 0.5)(1.85,0.4)(1.75,0.2)(1.73,.1)
%\psellipse[linestyle=dashed](4,0)(0.2,.65)
%\psellipse[](4,0)(0.2,.65)
\psellipticarc[](4,0)(0.2,.65){-90}{90}
\psellipticarc[linestyle=dashed](4,0)(0.2,.65){90}{270}
%
%\psccurve[](1.72,0.0)(1.73,0.50)(2.3,1.08)(2.27,0.50)(2.28,0.0)(2.27,-0.50)(1.7,-1.27)(1.73,-0.50)
%\psccurve[showpoints=true, linecolor=red](1.45,-1.4)(1.73,-0.50)(1.72,0.0)(1.73,0.50)(2.3,1.08)(2.27,0.50)(2.28,0.0)(2.27,-0.50)

\pscurve[linestyle=dashed,dash=.5pt 1pt]
(3.8,-0.7)(3.7,-.6)(3.69,-.5)(3.7,-.4)(3.75,-.3)(3.8,-.15)
(3.8,0)(3.85,.15)(4,.25)(4,.25)(4.15,.4)(4.2,.5)(4.3,.58)

\pscurve[](4.3,.6)(4.4,0.52)(4.2,0)(4.18,-.2)(4.15,-.3)(4,-.56)(3.8,-0.69)
%%%%%%%%%%%%%%%%%%%%%%%%%%%%%%%%%%%%%%%%%%%%%%%%%%%%%%%
%%%%%%%%%%%%%%%%%%%%%%%%%%%%%%%%%%%%%%%%%%%%%%%%%%%%%
%\psellipse[linecolor=red](1.58,-.36)(0.69,.6)
%\psellipticarc[linecolor=red](1.58,-.36)(1,.6){-100}{0}
%%%%%%%%%%%%%%%%%%%%%%%%%%%%%%%%%%%%%%%%%%%%%%%%%%%%%%%%
%%%%%%%%%%%%%%%%%%%%%%%%%%%%%%%%%%%%%%%%%%%%%%%%%%%%%%%%%
\uput{0.3}[0](0,2.8){$\cH$}
\uput{0.3}[0](1.6,1.8){$\cH^1_A=\Dom(A^{1/2})$}
\uput{0.3}[0](3.6,1.1){$\cH^2_A=\Dom(A)$}
\uput{0.3}[0](1.6,-1.8){$\cH^1_B=\Dom(|B|^{1/2})$}
\uput{0.3}[0](3.6,-1.1){$\cH^2_B=\Dom(B)$}
\end{pspicture}
\caption{\label{fig:2} The Sobolev-like scale of spaces for the operators $A$
and $B$. $\cH_A^2\ne\cH_B^2$ but the form-domain stability condition
$\cH_A^1=\cH_B^1$ still holds.}
\end{figure}

\begin{remark}
If $\Dom(A)\neq \Dom(B)$, then any of the possibilities $\Dom(B)\subset\Dom(A)$ and $\Dom(B)\triangle\Dom(A)\neq \emptyset$ may occur. Indeed, in the Hilbert space $\fH=\ell^2(\N;\C^2)$ consider the self-adjoint operator $A$ defined in \eqref{A:def}. Let
\begin{equation*}
H=\bigoplus_{k\in\N}\begin{pmatrix} 1 & 1 \\ 1 & -q\end{pmatrix}\quad\text{with}\quad q=0\quad\text{or}\quad q=1
\end{equation*}
and $B=A^{1/2} H A^{1/2}$ defined on its natural domain. It is straightforward to verify that
\begin{equation*}
\Dom(B) = \ell^{2,2}(\N)\oplus\ell^{2,4}(\N)\subsetneq\Dom(A)\quad\text{if}\quad q=1
\end{equation*}
and
\begin{equation*}
\Dom(B) = \ell^{2,2}(\N)\oplus\ell^{2,2}(\N)\not\subset\Dom(A)\quad\text{if}\quad q=0.
\end{equation*}
In the second case we, obviously, have $\Dom(B)\triangle\Dom(A)\neq \emptyset$.
\end{remark}

\vskip .5cm

Revisiting Example \ref{exam}, one can illustrate the statement of Theorem
\ref{lem:n} as follows.  By direct computations one easily checks that
\vskip .5cm
\begin{itemize}
\item[(a)] the sets  $\Dom(A^{1/2})$ and $\Dom(|B|^{1/2})$ are in general position, that is, the symmetric difference
$\Dom(A^{1/2})\triangle \Dom(|B|^{1/2})$ is a non-empty set and,
 hence, (i), (ii), and (ii$'$) do not hold;
\vskip .5cm
\item[(b)] the operators
\begin{equation*}
X=\bigoplus_{k\in\N} \begin{pmatrix} k & 0 \\ 0 & k^{-1} \end{pmatrix}\quad\text{and}\quad Y=\bigoplus_{k\in\N} \begin{pmatrix} k^{-1} & 0 \\ 0 & k \end{pmatrix}
\end{equation*}
are, obviously, unbounded. Hence (iii) and (iii$'$) do not hold;
\vskip .5cm
\item[(c)] the involution
\begin{equation*}
 K= \bigoplus_{k\in\N} \begin{pmatrix} 0 & k^{-1} \\ k & 0 \end{pmatrix}
\end{equation*}
is obviously unbounded so that (iv) does not hold;
\end{itemize}

\medskip

\noindent and, finally,

\medskip

\begin{itemize}
\item[(d)] $\sign(B) \Dom(A^{1/2})$ is not a subset of $\Dom(A^{1/2})$ and, hence, (v) fails to
hold.
\end{itemize}

\bigskip

It is also worth mentioning that the $A$-form boundedness of the self-adjoint operator  $B$ guaranteed
under Hypothesis \ref{hyp1} by Remark \ref{rem:1}
does not imply the $A$-form boundedness of $|B|$ in general (cf.~Theorem
\ref{lem:n} (iii)) which seems to be a bit unexpected.

\begin{figure}[htb]
\centering
\begin{pspicture}[](0,-3)(7,3)
\pscurve[linecolor=black,linewidth=1pt]
(0,2.5)
(2,1.2)
(4,0.65)
(7,.5)
%\psecurve[linestyle=dashed](3,4)(4,4.23)(5,4.43)
\pscurve[linecolor=black,linewidth=1pt]
(0,-2.5)
(2,-1.2)
(4,-0.65)
(7,-.5)
%\psecurve[linestyle=dashed](3,4)(4,4.23)(5,4.43)
\psellipse[](0,0)(0.5,2.5)
\psellipticarc[linestyle=dashed](2,0)(0.3,1.2){90}{270}
\psellipticarc[](2,0)(0.3,1.2){-90}{90}

\pscurve[linestyle=dashed,dash=.5pt 1pt](1.46,-1.48)(1.29,-1.44)(1.18,-1.28)(1.73,-0.50)(1.72,-.25)

%\pscurve[](1.48,-1.48)(1.7,-1.3)(2,-.95)(2.28,0.0)(2.27,0.50)(2.33,0.8)(2.43,.934)(2.65,0.93)
\pscurve[](1.48,-1.48)(1.7,-1.3)(2,-.95)(2.28,0.0)(2.27,0.50)(2.43,0.8)(2.43,.934)(2.40,1.01)

%\pscurve[linestyle=dashed,dash=.5pt 1pt](2.65,0.93)
%(2.7,.9)(2.72,.75)(2.71,.7)(2.7,0.65)(2.6,.55)(2.54, 0.52)(2.45, 0.51)(2.3, 0.5)(2.1, 0.5)(1.85,0.4)(1.75,0.2)(1.73,.1)

\pscurve[linestyle=dashed,dash=.5pt 1pt](2.40,1.01)(2.35, 0.99)(2.3, 0.98)(2.1, 0.8)(1.85,0.4)(1.75,0.2)(1.73,.1)

%\psellipse[linestyle=dashed](4,0)(0.2,.65)
%\psellipse[](4,0)(0.2,.65)
\psellipticarc[](4,0)(0.2,.65){-90}{90}
\psellipticarc[linestyle=dashed](4,0)(0.2,.65){90}{270}
%
%\psccurve[](1.72,0.0)(1.73,0.50)(2.3,1.08)(2.27,0.50)(2.28,0.0)(2.27,-0.50)(1.7,-1.27)(1.73,-0.50)
%\psccurve[showpoints=true, linecolor=red](1.45,-1.4)(1.73,-0.50)(1.72,0.0)(1.73,0.50)(2.3,1.08)(2.27,0.50)(2.28,0.0)(2.27,-0.50)

\pscurve[linestyle=dashed,dash=.5pt 1pt]
(3.8,-0.7)(3.7,-.6)(3.69,-.5)(3.7,-.4)(3.75,-.3)(3.8,-.15)
(3.8,0)(3.85,.15)(4,.25)(4,.25)(4.15,.4)(4.2,.5)(4.3,.58)

\pscurve[](4.3,.6)(4.4,0.52)(4.2,0)(4.18,-.2)(4.15,-.3)(4,-.56)(3.8,-0.69)
%%%%%%%%%%%%%%%%%%%%%%%%%%%%%%%%%%%%%%%%%%%%%%%%%%%%%%%
%%%%%%%%%%%%%%%%%%%%%%%%%%%%%%%%%%%%%%%%%%%%%%%%%%%%%
%\psellipse[linecolor=red](1.58,-.36)(0.69,.6)
%\psellipticarc[linecolor=red](1.58,-.36)(1,.6){-100}{0}
%%%%%%%%%%%%%%%%%%%%%%%%%%%%%%%%%%%%%%%%%%%%%%%%%%%%%%%%
%%%%%%%%%%%%%%%%%%%%%%%%%%%%%%%%%%%%%%%%%%%%%%%%%%%%%%%%%
\uput{0.3}[0](0,2.8){$\cH$}
\uput{0.3}[0](1.6,1.8){$\cH^1_A=\Dom(A^{1/2})$}
\uput{0.3}[0](3.6,1.1){$\cH^2_A=\Dom(A)$}
\uput{0.3}[0](1.6,-1.8){$\cH^1_B=\Dom(|B|^{1/2})$}
\uput{0.3}[0](3.6,-1.1){$\cH^2_B=\Dom(B)$}
\end{pspicture}
\caption{\label{fig:3} The Sobolev-like scale of spaces  for the
operators $A$ and $B$. The form-domain stability condition
$\cH_A^1=\cH_B^1$ does not hold and  the domains $\Dom(A^{1/2})$
and $\Dom(|B|^{1/2})$ are in general position in accordance with
Theorem   \ref{lem:n} (i), (ii)  and (ii$'$),
that is, $\cH_A^1\triangle\cH_B^1\ne \emptyset$.}
\end{figure}
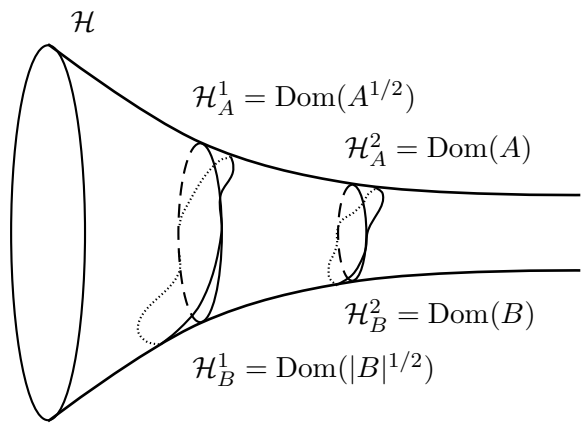

The corresponding  (typical) diagram illustrating ``counterexample'' \ref{exam}
is depicted in Fig.~\ref{fig:3}.

%%%%%%%%%%%%%%%%%%%%%%%%%%%%%%%%%%%%%%%%%%%%%%%%%%%%%%%%%%%%%%%%%%%%%%%%%%%
\section{On a Converse to the First Representation Theorem}\label{sec:conv}
%%%%%%%%%%%%%%%%%%%%%%%%%%%%%%%%%%%%%%%%%%%%%%%%%%%%%%%%%%%%%%%%%%%%%%%%%%%

In the semibounded case there is a one-to-one correspondence between the closed symmetric forms and the associated self-adjoint operators. For non-semibounded case the situation may be quite different and examples
of an operator associated with infinitely many sesquilinear forms naturally arise.

To illustrate this phenomenon we assume the following hypothesis.

\begin{hypothesis}\label{hypo:nonun}
Assume that  the (separable) Hilbert space $\fH$ admits an orthogonal decomposition
 $\fH=\fH_0\oplus\fH_1$ such that  the subspaces $\fH_0$ and $\fH_1$ have infinite dimension.
Suppose that  $D:\fH_1\to\fH_0$ is a closed densely defined operator.
Assume, in addition, that $D$ has a bounded inverse and let $D=U|D|$ be the polar decomposition
of the operator $D$ with a unitary $U:\fH_0\to\fH_1.$
\end{hypothesis}

Given $\mu\in[0,1]$, under Hypothesis \ref{hypo:nonun} introduce the self-adjoint positive definite operator matrix
\begin{equation*}
A_\mu=
\begin{pmatrix}
|D|^{2-2\mu} &0\\
0&|D^*|^{2\mu}
\end{pmatrix}
\end{equation*}
on
\begin{equation*}
\Dom(A_\mu)=\Dom(|D|^{2-2\mu})\oplus
\Dom(|D^*|^{2\mu})
\end{equation*}
and the self-adjoint bounded involution on $\fH_0\oplus \fH_1$
\begin{equation*}
H=\begin{pmatrix}
0 &U\\
U^*&0
\end{pmatrix}.
\end{equation*}

Clearly, the  self-adjoint operators  $A_\mu$ and $H$ satisfy  Hypothesis \ref{hyp1} and, therefore,
the sesquilinear symmetric form
\begin{equation}\label{bmu}
\fb_\mu[x,y]  := \langle A^{1/2}_\mu x, H A^{1/2}_\mu y \rangle_{\fH},
\quad x,y\in \Dom[\fb_\mu]=\Dom(A^{1/2}_\mu),
\end{equation}
is a $0$-closed form in the sense of McIntosh (cf.~Remark \ref{rem:mcintosh}). Since
\begin{equation*}
\Dom(A^{1/2}_\mu)\ne \Dom(A^{1/2}_\nu), \quad \mu\ne \nu,
\end{equation*}
the $0$-closed forms $\fb_\mu$ are defined on different domains
and, therefore, $\fb_\mu\ne \fb_\nu$ whenever $\mu$ and $\nu$, with
$\mu,\nu\in[0,1]$,
 are different.

By the First Representation Theorem, there exists
a unique self-adjoint operator $B_\mu$ associated with the form
$\fb_\mu$. However, our next result shows that this operator does not depend on $\mu$ and, therefore,
there exist infinitely many $0$-closed forms the self-adjoint operator
$B:=B_\mu$ is associated with.

\begin{proposition}\label{propo:nonun}
Assume Hypothesis \ref{hypo:nonun}. Then
\begin{itemize}
\item[(i)] The block operator matrix
\begin{equation}\label{B:def}
B:=\begin{pmatrix} 0 & D \\ D^\ast & 0 \end{pmatrix}
\end{equation}
defined on its natural domain $\Dom(B)=\Dom(D^*)\oplus\Dom (D)$ is a self-adjoint operator.
\item [(ii)] For any $\mu\in [0,1]$ the  operator $B$  is associated
with the  form $\fb_\mu$ given by \eqref{bmu}.
\item [(iii)]
 The form $\fb_\mu$ is represented by the operator $B$ if and only if $\mu=1/2$.
\end{itemize}
\end{proposition}

\begin{proof}
Under Hypothesis \ref{hypo:nonun} let
\begin{equation*}
D=U|D|\quad \text{on } \Dom(D)=\Dom(|D|)
\end{equation*}
be the polar decomposition of $D$ (cf.~\cite[Sect.~VI.2.7]{Kato}). Recall that
\begin{equation*}
D^*=U^*|D^*|\quad \text{on } \Dom(D^*)=\Dom(|D^*|).
\end{equation*}
By a result in \cite[Theorem 2.7]{Gesztesy}, for any $\mu\in[0,1]$
the operators $D$ and $D^*$ can be represented as the products
\begin{equation}\label{D}
D = |D^\ast|^\mu U |D|^{1-\mu}\quad\text{on}\quad \Dom(D)
\end{equation}
and
\begin{equation}\label{D*}
D^* = |D|^{1-\mu} U^* |D^*|^{\mu}\quad\text{on}\quad \Dom(D^*) .
\end{equation}
Therefore, the operator matrix
\begin{equation*}
B=
\begin{pmatrix}
0&D\\
D^*&0
\end{pmatrix}\quad \text{ on } \Dom (B)=\Dom(D^*)\oplus \Dom(D)
\end{equation*}
admits the factorization
\begin{equation*}
B=
\begin{pmatrix}
|D|^{1-\mu} &0\\
0&|D^*|^\mu
\end{pmatrix}
\begin{pmatrix}
0 &U\\
U^*&0
\end{pmatrix}\begin{pmatrix}
|D|^{1-\mu} &0\\
0&|D^*|^\mu
\end{pmatrix}.
\end{equation*}

By Lemma \ref{lem:1} the operator
\begin{equation*}
B_\mu=A_\mu^{1/2}HA_\mu^{1/2}
\end{equation*}
defined on its natural maximal domain
\begin{equation*}
\Dom(B_\mu)=\{x\in\Dom(A^{1/2}_\mu)\,|\, HA^{1/2}_\mu x\in\Dom(A^{1/2}_\mu)\}
\end{equation*}
is a self-adjoint operator with a bounded inverse.

It remains to observe that due to \eqref{D} and \eqref{D*},  $\Dom(B)=\Dom(B_\mu)$ and, therefore,
\begin{equation*}
B=B_\mu\quad \text{for all }\mu\in [0,1].
\end{equation*}

By Theorem \ref{thm:main}, the self-adjoint operator $B_\mu$ is associated with the form $\fb_\mu$, so does $B$ which proves (i) and (ii).

(iii). The claim follows from  Theorem \ref{lem:n} (i) and Theorem \ref{thm:main:2}.
\end{proof}

We illustrate the statement of Proposition \ref{propo:nonun} on the classical example of the
free Dirac operator of Quantum Mechanics.

\begin{example}\label{Dirac}
Let $\fH_0 \cong \fH_1 = L^2(\R^3; \C^2)$.
Consider the free Dirac operator defined on its natural domain as a block operator matrix
\begin{equation*}
B:=\begin{pmatrix} I & \partial  \\ \partial^\ast & - I \end{pmatrix},
\end{equation*}
where $\partial =\ii \vec{\sigma}\cdot \vec{\nabla}$ and
$\vec{\sigma}=(\sigma_1, \sigma_2,\sigma_3)$, the Pauli matrices,
\begin{equation*}    \sigma_1 = \begin{pmatrix} 0&1\\ 1&0 \end{pmatrix},
   \quad  \sigma_2 =\begin{pmatrix} 0&-\ii\\ \ii&0 \end{pmatrix},
    \quad \sigma_3 =  \begin{pmatrix} 1&0\\ 0&-1 \end{pmatrix}.
\end{equation*}

We emphasize that the ``pathology'' we dealt with in Proposition \ref{propo:nonun},
already occurs for the Dirac operator. Indeed, it is well known that the Dirac operator $B$ is self-adjoint on its natural domain and that $B$ has absolutely continuous spectrum of infinite multiplicity
filling in the set $(-\infty, -1]\cup [1,\infty).$ By the Spectral Theorem,
the operator $B$ is unitarily equivalent to the block operator matrix
\begin{equation*}
\widetilde B=\begin{pmatrix} 0&M\\
M^*&0\end{pmatrix},
\end{equation*}
where $M$ is the multiplication operator by the independent variable
in the Hilbert space
\begin{equation*}
\fL=L^2((1,\infty);\fH^\prime)
\end{equation*}
of vector-valued functions with values in an infinite dimensional (separable)
Hilbert space $\fH^\prime$.

By Proposition \ref{propo:nonun}, the operator $\widetilde B$ is
associated with infinitely many $0$-closed sesquilinear symmetric forms, so does the Dirac operator $B$.
\end{example}

%%%%%%%%%%%%%%%%%%%%%%%%%%%%%%%%%%%%%%%%%%%%%%%%

\end{document}